\documentclass{amsart}
\usepackage[dvips]{graphicx}
\usepackage{amscd}
\usepackage{amsmath}
\usepackage{amsxtra}
\usepackage{amsfonts}
\usepackage{amssymb}
\newtheorem{theorem}{Theorem}[section]
\newtheorem{corollary}[theorem]{Corollary}

\newtheorem{proposition}[theorem]{Proposition}
\theoremstyle{definition}
\newtheorem{definition}[theorem]{Definition}
\newtheorem{remark}[theorem]{Remark}

\newtheorem{example}[theorem]{Example}
\theoremstyle{remark}

\renewcommand{\theclaim}{\textup{\theclaim}}

\newtheorem*{acknowledgements}{Acknowledgements}

\numberwithin{equation}{section}

\def\openone

{\mathchoice

{\hbox{\upshape \small1\kern-3.3pt\normalsize1}}

{\hbox{\upshape \small1\kern-3.3pt\normalsize1}}

{\hbox{\upshape \tiny1\kern-2.3pt\SMALL1}}

{\hbox{\upshape \Tiny1\kern-2pt\tiny1}}}

\makeatletter

\newbox\ipbox

\newcommand{\ip}[2]{\left\langle #1\,|\,#2\right\rangle}

\newcommand{\diracb}[1]{\left\langle #1\mathrel{\mathchoice

{\setbox\ipbox=\hbox{$\displaystyle \left\langle\mathstrut
#1\right.$}

\vrule height\ht\ipbox width0.25pt depth\dp\ipbox}

{\setbox\ipbox=\hbox{$\textstyle \left\langle\mathstrut #1\right.$}

\vrule height\ht\ipbox width0.25pt depth\dp\ipbox}

{\setbox\ipbox=\hbox{$\scriptstyle \left\langle\mathstrut
#1\right.$}

\vrule height\ht\ipbox width0.25pt depth\dp\ipbox}

{\setbox\ipbox=\hbox{$\scriptscriptstyle \left\langle\mathstrut
#1\right.$}

\vrule height\ht\ipbox width0.25pt depth\dp\ipbox}

}\right. }

\newcommand{\dirack}[1]{\left. \mathrel{\mathchoice

{\setbox\ipbox=\hbox{$\displaystyle \left.\mathstrut
#1\right\rangle$}

\vrule height\ht\ipbox width0.25pt depth\dp\ipbox}

{\setbox\ipbox=\hbox{$\textstyle \left.\mathstrut #1\right\rangle$}

\vrule height\ht\ipbox width0.25pt depth\dp\ipbox}

{\setbox\ipbox=\hbox{$\scriptstyle \left.\mathstrut
#1\right\rangle$}

\vrule height\ht\ipbox width0.25pt depth\dp\ipbox}

{\setbox\ipbox=\hbox{$\scriptscriptstyle \left.\mathstrut
#1\right\rangle$}

\vrule height\ht\ipbox width0.25pt depth\dp\ipbox}

} #1\right\rangle}

\newcommand{\lonet}{L^{1}\left(  \mathbb{T}\right)}

\newcommand{\linft}{L^{\infty}\left(  \mathbb{T}\right)}

\newcommand{\loner}{L^{1}\left(\mathbb{R}\right)}

\newcommand{\ltwor}{L^{2}\left(\mathbb{R}\right)}

\newcommand{\cj}[1]{\overline{#1}}

\input cyracc.def

\begin{document}
\title[MRA super-wavelets]{MRA super-wavelets}
\author{Stefan Bildea, Dorin Ervin Dutkay, Gabriel Picioroaga}
\address{Department of Mathematics\\
The University of Iowa\\
14 MacLean Hall\\
Iowa City, IA 52242-1419\\
U.S.A.\\} \email{Stefan Bildea: sbildea@math.uiowa.edu} \email{Dorin
Ervin Dutkay: ddutkay@math.rutgers.edu} \email{Gabriel Picioroaga:
gpicioro@math.uiowa.edu}
\thanks{}
\subjclass{42C40, 37C30, 42C30} \keywords{Multiresolution,
wavelet, low-pass filter, scaling function, transfer operator,
cascade algorithm, representation}

\begin{abstract}
We construct a multiresolution theory for
$\ltwor\oplus...\oplus\ltwor$. For a good choice of the dilation
and translation operators on these larger spaces, it is possible
to build singly generated wavelet bases, thus obtaining
multiresolution super-wavelets. We give a characterization of
super-scaling function, we analyze the convergence of the cascade
algorithms and give examples of super-wavelets. Our analysis
provides also more insight into the Cohen and Lawton condition
for the orthogonality of the scaling function in the classical
case on $\ltwor$.
\end{abstract}

\maketitle \tableofcontents
\section{\label{Intro}Introduction}
\par
The applications of wavelet theory to signal processing and image
processing are now well known. Probably the main reason for the
success of the wavelet theory was the introduction of the concept of
multiresolution analysis (MRA), which provided the right framework
to construct orthogonal wavelet bases with good localization
properties.
\par
One of the problems in networking is {\it multiplexing}, which consists of sending multiple signals or streams of information on a carrier at the same time in the form of a single, complex signal and then recovering the separate signals at the receiving end.
\par
In \cite{HL}, Deguang Han and David Larson have shown that the
technique of multiresolution analysis breaks down, when multiplexing
is required, if one just amplifies (in the operator theory sense)
the steps used in the construction of MRA-wavelets. We will state
this more precisely in a moment.
\par
In this paper we will show how multiresolution constructions {\em
can} be realized for multiple signals, provided some slight
modifications are done to the usual dilation and translation
operators. We believe that our constructions have potential for
applications to multiplexing problems.
\par
In \cite{Dut2} and \cite{Dut3}, the second named author introduced a
certain affine structure on the space $\ltwor\oplus...\oplus\ltwor$
which was shown to admit multiresolution wavelet bases. In this
paper, we will further analyze this affine structure, give
characterizations of its scaling vectors (Theorem \ref{th1_9},
Corollary \ref{cor1_15}). Then, the well known conditions for the
orthogonality of the scaling functions due to A. Cohen \cite{Co90}
and W. Lawton \cite{Law91a} are extended to this larger space in
Theorem \ref{th2_1}. We study the convergence of the cascade
operator which provides numerical approximations for the scaling
function. Finally, in Section 5, we construct several examples of
super-scaling functions and super-wavelets and we prove that all
these structures admit super-wavelet bases.
\par
In this introduction we recall several fundamental ideas and
notions of the wavelet theory. We refer to \cite{Dau92},
\cite{BraJo} for more information on the topic.
\par  A {\it
wavelet} is a function $\psi\in\ltwor$ with the property that
$$\{U^mT^n\psi\,|\, m,n\in\mathbb{Z}\}$$
is an orthonormal basis for $\ltwor$, where $U$ is the dilation
operator
$$Uf(x)=\frac{1}{\sqrt{2}}f\left(\frac{x}{2}\right),\quad(f\in\ltwor,x\in\mathbb{R}),$$
and $T$ is the translation operator
$$Tf(x)=f(x-1),\quad(f\in\ltwor,x\in\mathbb{R}).$$
\par
A {\it multiresolution analysis} is an increasing nest of closed
subspaces $(V_n)_{n\in\mathbb{Z}}$ of $\ltwor$ which has a dense
union, a zero intersection, $UV_n=V_{n-1}$, and there is a vector
$\varphi$ in $V_0$ such that
$$\{T^k\varphi\,|\, k\in\mathbb{Z}\}$$
is an orthonormal basis for $V_0$. This $\varphi$ is then called
an {\it orthogonal scaling function}.
\par
The wavelets $\psi$ are then constructed from the scaling function
in such a way that
$$\{T^k\psi\,|\, k\in\mathbb{Z}\}$$
is an orthonormal basis for $W_0:=V_1\ominus V_0$, the relative
orthocomplement.
\par
There are generalizations of the multiresolution analyses some of
which we will turn to later. One of the directions for these
generalizations (see \cite{BCMO}) is to replace the Hilbert space
$\ltwor$ by an abstract one $H$ and the dilation and translation
operators $U$ and $T$ by some abstract unitaries that verify a
commutation relation such as
$$UTU^{-1}=T^2.$$
The two unitaries will generate what we call a wavelet
representation. We adopt this representation theoretic view point,
and our approach emphasizes the strong connections between wavelet
theory and the spectral properties of a certain transfer operator
(Definition \ref{def0_0}).
\par
We maintain an abstract flavor throughout the paper, but we concentrate on some concrete examples
of wavelet representations on the
Hilbert space $\ltwor\oplus...\oplus\ltwor$ (finite sum), which we
will describe in a moment.
\par
It is known (see \cite{HL}, Proposition 5.16), that no orthogonal
scaling function can be constructed for $\ltwor\oplus\ltwor$ with
the dilation $U\oplus U$ and the translation $T\oplus T$, therefore
multiplexing can not be obtained in this way. However, we will see
that a multiresolution theory can be developed on
$\ltwor\oplus...\oplus\ltwor$ with plenty of examples of orthogonal
scaling functions and wavelets, if some slight modifications are
done to these operators.
\par
The dilation and translation operators are constructed as follows:
first consider a fixed cycle, that is, a periodic orbit for the
map $z\mapsto z^2$ on the unit circle. Let
$C:=\{z_1,z_2,...,z_p\}$ be this cycle,
$z_1^2=z_2,z_2^2=z_3,...,z_{p-1}^2=z_p,z_p^2=z_1$. The Hilbert
space is
$$H_C:=\underbrace{\ltwor\oplus...\oplus\ltwor}_{p\mbox{ times
}},$$ the translation operator $T_C$ on $H_C$ is given by
$$T_C(\xi_1,...,\xi_p)=(z_1T\xi_1,...,z_pT\xi_p),\quad(\xi_1,...,\xi_p\in\ltwor),$$
the dilation operator $U_C$ on $H_C$ is defined by
$$U_C(\xi_1,...,\xi_p)=(U\xi_2,U\xi_3,...,U\xi_p,U\xi_1),\quad(\xi_1,...,\xi_p\in\ltwor),$$
where $T$ and $U$ are the translation and dilation operators on
$\ltwor$ defined before.
\par
It was shown in \cite[Proposition 2.13]{Dut2} that if a filter
$m_0\in\linft$ satisfies the condition
$$|m_0(z_i)|=\sqrt{N},\quad(i\in\{1,...,p\}),$$ for some cycle $C=\{z_1,...,z_p\}$ then there is a scaling vector $\varphi\in H_C$.
Such a cycle is called an {\it $m_0$-cycle}. \par With the operators
$U_C,T_C$, multiresolutions can be constructed and super-wavelets
(meaning wavelets in spaces larger then $\ltwor$) are obtained. The
representations that correspond to different cycles are disjoint so
that we will see that scaling functions and multiresolutions can be
constructed also for the direct sum of these representations.
Theorem 5.3 shows that, gathering {\it all} $m_0$-cycles, one
obtains an orthogonal scaling vector in the orthogonal sum $\oplus_C
H_C$.
\par
We prove here that the scaling vector constructed out of some
$m_0$-cycles are orthogonal only when all the $m_0$-cycles are taken
into consideration, and when this is not the case, the
super-wavelets form a normalized tight frame (Theorem \ref{th2_1}).
In particular, one of the consequences of our analysis is that the
MRA normalized tight frame wavelets, which are known to occur on
$\ltwor$ when the low-pass filter is not judiciously chosen, are in
fact projections of good, orthogonal MRA super-wavelets. This is due
to the fact that the only cycle considered in the construction of
wavelets on $\ltwor$ is the trivial cycle $\{1\}$, while the filter
$m_0$ might have other cycles. The direct relation to Cohen's
orthogonality condition \cite{Co90} is now clear.
\par
Our paper is structured as follows:\par In Section \ref{wave} we
define the term of wavelet representation (Definition \ref{def1_1})
and present the main examples that we will work with (especially
Examples \ref{ex1_5} and \ref{ex1_6}).
\par
In Section \ref{scal} we define the notions of multiresolution
analysis and scaling vector (Definition \ref{def1_6}), imitating the
one existent for $\ltwor$ and give a characterization theorem for
scaling vectors (Theorem \ref{th1_9}). The theorem is then applied
to our main examples, and we obtain in this way conditions which
characterize scaling vectors in $\ltwor\oplus...\oplus\ltwor$
(Corollary \ref{cor1_15}).
\par
Once a scaling vector is found, the construction of wavelets can be
done in exactly the same fashion as for $\ltwor$ (see Proposition
\ref{prop1_11}).
\par
In Section 3 we will see how one can obtain {\it super-scaling
functions} and{\it  super-wavelets} from a trigonometric polynomial
filter. Each component of the super-scaling function will be
compactly supported. Just like in the $\ltwor$ case, the
super-wavelet usually generates a normalized tight frame and, to get
an orthogonal basis, some extra conditions must be imposed on the
initial low-pass filter from which the super-wavelet is constructed
(such as the Cohen condition or Lawton's condition, see Theorem
\ref{th2_1}).
\par
The scaling vector is approximated by the so-called {\it cascade
algorithm}. One starts with a well chosen function and then the {\it
cascade operator} is applied successively to it. In this way one
obtains a sequence which approaches the scaling vector in norm. We
will see in Section \ref{casc} how the initial function must be
chosen so that the algorithm is convergent.
\par
Abstract or geometric constructions in wavelet theory must be tested
against examples and explicit algorithms. We end our paper with
Section \ref{exam} which contains several examples showing that
plenty of multiresolution super-wavelets can be constructed.
\par
Some notations that we will use in this paper:
\\The Fourier transform of $f\in\loner$ is given by
$$\widehat{f}(\xi)=\int_{\mathbb{R}}f(x)e^{-i\xi x}\,dx.$$
We denote by $\mathbb{T}$ the unit circle $\{z\in\mathbb{C}\,|\,
|z|=1\}$ and by $\mu$, the normalized Haar measure on
$\mathbb{T}$. We often identify functions $f$ on $T$ with
$2\pi$-periodic functions on $\mathbb{R}$ or with functions on the
interval $[-\pi,\pi]$. The identification is given by
$$f(z)\leftrightarrow f(\theta)\mbox{ where } z=e^{-i\theta}.$$
\par
Many key properties of the scaling vectors are encoded in the
spectral properties of a certain transfer operator. This is
defined as follows:
\par
\begin{definition}\label{def0_0}
Let $N\geq2$ be an integer. The {\it transfer operator} is
associated to a function $m_0\in\linft$ and is defined by
$$R_{m_0}f(z)=\frac{1}{N}\sum_{w^N=z}|m_0(w)|^2f(w),\quad(z\in\mathbb{T},f\in\lonet).$$
\par
A function $h\in\lonet$ is called {\it harmonic} with respect to
$R_{m_0}$ if
$$R_{m_0}h=h.$$
\end{definition}

\section{\label{wave}Wavelet representations}
\par Fix an integer $N\geq2$ called the {\it scale}.
 Wavelet representations are an abstract version of the situation
 existent on $\ltwor$ as explained in the introduction. The
 Hilbert space $\ltwor$ is replaced by an abstract one $H$ and the
 dilation and translation operators are replaced by two unitaries
 $U$ and $T$ satisfying the following commutation relation
 $$UTU^{-1}=T^N.$$
 The translation generates a representation of $\linft$ by the
 Borel functional calculus. We give here the proper definition of
a wavelet representation and present some examples.
\begin{definition}\label{def1_1}
A {\it wavelet representation} is a triple
$\tilde{\pi}:=(H,U,\pi)$ where $H$ is a Hilbert space, $U$ is a
unitary on $H$ and $\pi$ is a representation of $\linft$ on $H$
such that
\begin{equation}\label{eqcov}
U\pi(f)U^{-1}=\pi(f(z^N)),\quad(f\in\linft). \end{equation}
(here, by $f(z^N)$ we mean the map $z\mapsto f(z^N)$).
\par
A wavelet representation is called {\it normal} if for any sequence
$(f_n)_{n\in\mathbb{N}}$ which converges pointwise a.e. to a
function $f\in\linft$ and such that $\|f_n\|_{\infty}\leq M$,
$n\in\mathbb{N}$ for some $M>0$, the sequence $\{\pi(f_n)\}$
converges to $\pi(f)$ in the strong operator topology.
\par
Sometimes we call $U$ the dilation and $T:=\pi(z)$ the translation
of the wavelet representation (here $z$ indicates the identity
function on $\mathbb{T}$, $z\mapsto z$).
\end{definition}

\begin{example}\label{ex1_2}
The main example of a wavelet representation is the classical one:
$H=\ltwor$,
$$U\xi(x)=\frac{1}{\sqrt{N}}\xi\left(\frac{x}{N}\right),\quad(\xi\in\ltwor),$$
and $\pi$ is defined by its Fourier transform
$${\widehat\pi}(f)(\xi)=f\xi,\quad(f\in\linft,\xi\in\ltwor).$$
In particular
$$\widehat{T}\xi(x)=e^{-ix}\xi(x),\mbox{ so
}T(\xi)(x)=\xi(x-1),\quad(\xi\in\ltwor,x\in\mathbb{R}).$$ We
denote this normal wavelet representation by $\mathfrak{R}_0$.
\end{example}

\begin{example}\label{def1_3}
If $(H_i,U_i,\pi_i)$ are (normal) wavelet representations for
$i\in\{1,...,n\}$, then
$(\oplus_{i=1}^nH_i,\oplus_{i=1}^nU_i,\oplus_{i=1}^n\pi_i)$ is a
(normal) wavelet representation called the direct sum of the given
wavelet representations.
\end{example}

\begin{example}\label{ex1_4}
We call {\it cycle} a set $\{z_1,...,z_p\}$ of distinct points in
$\mathbb{T}$, such that $z_1^N=z_2,z_2^N=z_3,...,z_p^N=z_1$. $p$
is called the length of the cycle. $\{1\}$ is called the trivial
cycle.
\par
Let $(H,U,\pi)$ be a (normal) wavelet representation. Let
$C:=\{z_1,...,z_p\}$ be a cycle and
$\alpha_1,...,\alpha_p\in\mathbb{T}$. Define
$$H_{C,\alpha}:=\underbrace{H\oplus H\oplus...\oplus H}_{p\mbox{ times}}$$
and, for $f\in\linft$, $\xi_1,...,\xi_p\in H$,
$$U_{C,\alpha}(\xi_1,...,\xi_p):=(\alpha_1U\xi_2,\alpha_2
U\xi_3,...,\alpha_{p-1}U\xi_{p},\alpha_pU\xi_1),,$$
$$\pi_{C,\alpha}(f)(\xi_1,...,\xi_p)=(\pi(f(z_1z))\xi_1,\pi(f(z_2z))\xi_2,...,\pi(f(z_pz))\xi_p).$$
Then
$\tilde{\pi}_{C,\alpha}:=(H_{C,\alpha},U_{C,\alpha},\pi_{C,\alpha})$
is a (normal) wavelet representation which we call the {\it cyclic
amplification} of $\tilde{\pi}$ with cycle $C$ and modulation
$\alpha$. We leave it to the reader to check that this is indeed a
(normal) wavelet representation (see also \cite{Dut2}).
\par
Note that the cyclic amplification with the trivial cycle and
$\alpha_1=1$ is the initial wavelet representation.
\par
When $\alpha_1=...=\alpha_p=1$ we will use also the notation
$\tilde{\pi}_C:=\tilde{\pi}_{C,\alpha}$.
\end{example}

\begin{example}\label{ex1_5}
If $C$ is a cycle of length $p$ and $\alpha_1,...,\alpha_p$ are in
$\mathbb{T}$, we denote by
$$\mathfrak{R}_{C,\alpha}=(\ltwor_{C,\alpha},U_{C,\alpha},\pi_{C,\alpha})$$
the cyclic amplification $(\mathfrak{R}_0)_{C,\alpha}$ of the main
representation $\mathfrak{R}_0$. When all $\alpha_i$ are $1$, we
use the shorter notation
$$\mathfrak{R}_C=(\ltwor_C,U_C,\pi_C).$$
\end{example}

\begin{example}\label{ex1_6}
The wavelet representation which is of main importance to us in
this paper is the direct sum of wavelet representations associated
to several cycles. That is, if $C_1,...,C_n$ are distinct cycles
and $\alpha_1,...,\alpha_n$ are some finite sets of numbers in
$\mathbb{T}$, then let $C:=C_1\cup...\cup C_n$,
$\alpha=(\alpha_1,...,\alpha_n)$ and define
$$\mathfrak{R}_{C,\alpha}:=\mathfrak{R}_{C_1,\alpha_1}\oplus...\oplus\mathfrak{R}_{C_n,\alpha_n},$$
which we will call the wavelet representation associated to the
cycles $C_1,...,C_n$ and the numbers $\{\alpha_i\}$.
\end{example}

\section{\label{scal}Scaling vectors and compactly supported super-wavelets}
We define now the key concepts of a multiresolution analysis and
scaling vectors. The definition generalize the existent ones on
$\ltwor$ (see \cite{Dau92}).
\begin{definition}\label{def1_6} Let $\left(H,U,\pi\right)$ be a wavelet representation.
A sequence $\left(V_n\right)_{n\in\mathbb{Z}}$ of closed subspaces
of $H$ with the properties
\begin{enumerate}
\item $V_n\subset V_{n+1}$;
\item $\overline{\cup_{n\in\mathbb{Z}} V_n}=H$;
\item $\cap_{n\in\mathbb{Z}} V_n=\{0\}$;
\item $U(V_n)=V_{n-1}$;
\item There is a $\varphi\in V_0$ such that
$\{T^k\varphi\,|\,k\in\mathbb{Z}\}$ is an orthonormal basis for
$V_0$,
\end{enumerate}
is called a {\it multiresolution analysis (MRA)}.
\par
 A vector $\varphi\in H$ for which there exists a MRA such that
(i)-(v) hold, (and $\varphi$ is as in (v)), is called {\it
orthogonal scaling vector}.
\end{definition}

\par
Before we give a characterization for orthogonal scaling vectors,
we need the following proposition:

\begin{definition}\label{prop1_7}
Let $\tilde{\pi}=(H,U,\pi)$ be a normal wavelet representation and
$v_1,v_2\in H$. The Radon-Nikodym derivative $h_{v_1,v_2}$ of the linear functional
$$f\mapsto\ip{\pi(f)v_1}{v_2}\quad(f\in\linft),$$
with respect to the Haar measure $\mu$ on $\mathbb{T}$, is called
the {\it correlation function} of $v_1$ and $v_2$. We also use
the notation $h_v:=h_{v,v}$. We can rewrite this as
$$\ip{\pi(f)v_1}{v_2}=\int_{\mathbb{T}}f\,h_{v_1,v_2}\,d\mu,\quad(f\in\linft).$$
\end{definition}
\par
The existence of the correlation function is guaranteed by the fact that the representation $\pi$ is normal.
\par
 If $\mathcal{S}$ is a
set of operators on a Hilbert space (such as $U$ plus $\pi(f)$ in
the case of a wavelet representation), then we denote by
$\mathcal{S}'$ its commutant (i.e. the set of all operators that
commute with all operators in $\mathcal{S}$).
\par
We proceed towards the theorem that gives a characterization of
orthogonal scaling vectors. The scaling vector will satisfy a
scaling equation ((ii) in the next theorem), which relates the
scaling vector to a {\it low-pass filter} $m_0\in\linft$. In
order to obtain a MRA, some non-degeneracy conditions must be
imposed on $m_0$. When $m_0$ {\it is} degenerate, a residual
subspace appears as the intersection of the multiresolution
subspaces $V_n$ (see \cite{BraJo97}).
\begin{definition}\label{def1_8}
A function $m_0\in\linft$ is called {\it degenerate} if
$|m_0(z)|=1$ for a.e. $z\in\mathbb{T}$ and there exists a
measurable function $\xi:\mathbb{T}\rightarrow\mathbb{T}$ and a
$\lambda\in\mathbb{T}$ such that
$$m_0(z)\xi(z^N)=\lambda\xi(z),\quad (z\in\mathbb{T}).$$
\end{definition}

\begin{theorem}\label{th1_9}
Let $\tilde{\pi}=\left(H,U,\pi\right)$ be a wavelet representation
and $\varphi\in H$. Then $\varphi$ is an orthogonal scaling vector
if and only if the following conditions are satisfied:
\begin{enumerate}
\item {\bf [Orthogonality]} The correlation function of $\varphi$
is $h_\varphi=1$ a.e.; \item {\bf [Scaling equation]} There
exists a non-degenerate $m_0\in\linft$ such that
$U\varphi=\pi(m_0)\varphi$; \item {\bf [Cyclicity]} There is no
proper projection $p$ in $\tilde{\pi}'$ such that
$p\varphi=\varphi$.
\end{enumerate}
\end{theorem}
\begin{proof} The ideas of the proof are similar to the case of
$\ltwor$, so we will only sketch them and refer also to
\cite{Dau92}, \cite{HeWe} and \cite{BraJo} for more details. \par
The orthogonality of the translates is equivalent to the fact
that the correlation function of $\varphi$ has Fourier
coefficients $\delta_0$, i.e., $h_\varphi$ is constant $1$. The
scaling equation is equivalent to the fact that $UV_0\subset
V_0$; the non-degeneracy of $m_0$ is equivalent to $\cap
V_n=\{0\}$ (see \cite{BraJo97} and \cite[Theorem 5.6]{Jor01}).
$\cup V_n$ is dense if and only if $\varphi$ is cyclic for
$\{U,\pi\}$ which is in turn equivalent to condition (iii)

\end{proof}

\begin{definition}\label{def1_12} Let $\left(H, \pi, U\right)$ be a wavelet representation. A vector $\varphi\in H$
is called a {\it scaling vector} with filter $m_0$ if it satisfies
conditions (ii) and (iii) in Theorem \ref{th1_9}, but here we allow
$m_0$ to be degenerate.
\end{definition}
\begin{remark}\label{rem1_12}
Note that the condition (iii) is equivalent to
$$\{U^{-m}\pi(f)\varphi\,|\,m\geq0\,f\in\linft\}\mbox{ is
dense in }H,$$ since the scaling equation (ii) in \ref{th1_9} and
the covariance relation (\ref{eqcov}) are satisfied.
\end{remark}
\par
We will apply the characterization theorem to the instance described
in Example \ref{ex1_6}. When applied to the classical wavelet
representation on $\ltwor$ we obtain a theorem similar to the one in
\cite{HeWe}, chapter 7.

\begin{corollary}\label{cor1_15}
Let $\mathfrak{R}_{C,\alpha}$ be the wavelet representation in
Example \ref{ex1_6}. Denote by $e^{-i\theta_{i,j}}$ the $j$-th point
of the cycle $C_i$. Then
$\varphi_{C}=\varphi_{C_1}\oplus\varphi_{C_2}\oplus...\oplus\varphi_{C_n}$
is an orthogonal scaling function for this wavelet representation if
and only if the following conditions are satisfied:
\begin{enumerate}
\item
$$ \sum_{i=1}^n\sum_{k=1}^{p_i}(\operatorname*{Per}(|\widehat
\varphi_{C_i,k}|^2))(\xi-\theta_{i,k})=1 \quad\mbox{for a.e.
}\xi\in\mathbb{R};$$ \item There exists a function $m_0\in\linft$
such that for a.e. $\xi\in\mathbb{R}$ and for all
$i\in\{1,...,n\}$:

$$ \alpha_{i,1} \sqrt{N}
\widehat\varphi_{C_i,2}(N\xi)=m_0(\theta_{i,1}+\xi)\widehat\varphi_{C_i,1}(\xi),$$
$$ \alpha_{i,2} \sqrt{N}
\widehat\varphi_{C_i,3}(N\xi)=m_0(\theta_{i,2}+\xi)\widehat\varphi_{C_i,2}(\xi),$$
$$...,$$
$$ \alpha_{i,p_i} \sqrt{N}
\widehat\varphi_{C_i,1}(N\xi)=m_0(\theta_{i,p_i}+\xi)\widehat\varphi_{C_i,p_i}(\xi);$$
\item For each $i\in\{1,...,n\},j\in\{1,...,p_i\}$, $\widehat\varphi_{C_i}$ does not vanish
on any subset $E$ of $\mathbb{R}$ invariant under dilations by
$N^{p_i}$ (i.e. $N^{p_i}E=E$) of positive measure.
\end{enumerate}
\end{corollary}
\begin{proof}
The formula in (i) is just the correlaton function for
$h_\varphi$. Applying the Fourier transform to the scaling
equation, one obtains (ii). The projections in the commutant of
this representation are given by sets which are invariant under
multiplication by $N^{p_i}$ (see \cite{Dut2}, lemma 2.14). The
fact that $m_0$ is non-degenerate is automatic. Indeed, suppose
not, then $|m_0(z)|=1$ for a.e $z\in\mathbb{T}$ so, for all
$i\in\{1,...,n\}$,
$$\sqrt{N^{p_i}}|\widehat{\varphi}_{C_i,1}(N^{p_i}\xi)|=|\widehat{\varphi}_{C_i,1}(\xi)|,\quad(\xi\in\mathbb{R}).$$
We will conclude that $\widehat{\varphi}_{C_i,1}$ must be $0$. If
for some $a>0$ there is a subset $E$ of
$[-N^{p_i},-1]\cup[1,N^{p_i}]$ of positive measure such that
$|\widehat{\varphi}_{C_i,1}(\xi)|\geq a$ for $\xi\in E$, then
$$|\widehat{\varphi}_{C_i,1}(\xi)|\geq\frac{a}{\sqrt{N^{kp_i}}},\quad\mbox{for
}\xi\in N^{kp_i}E,k\in\mathbb{Z}$$ But then this contradicts the
integrability of $\widehat{\varphi}_{C_i,1}$:
$$\int_{\mathbb{R}}|\widehat{\varphi}_{C_i,1}(\xi)|^2\,d\xi\geq\sum_{k\in\mathbb{Z}}\frac{a^2}{N^{kp_i}}N^{kp_i}\lambda(E)=\infty.$$
\end{proof}

\par
As soon as a scaling function is given for a wavelet
representation, the construction of wavelets follows the
procedure described for the $\ltwor$-case in \cite{Dau92}. We
present below the required ingredients in an abstract version.
For a proof look in \cite[Proposition 5.1]{Dut3}.

\begin{proposition}\label{prop1_11}
Let $\tilde{\pi}$ be a normal wavelet representation having an
orthogonal scaling function $\varphi$ with non-degenerate filter
$m_0$. Denote by $(V_n)_{n\in\mathbb{Z}}$ the associated MRA.
Assume that there are given the "high-pass filters"
$m_1,..,m_{N-1}\in\linft$ satisfying
\begin{equation}\label{eq1_2}
\frac{1}{\sqrt{N}}\left(\begin{array}{cccc}
m_0(z)&m_0(\rho z)&\hdots&m_0(\rho^{N-1}z)\\
m_1(z)&m_1(\rho z)&\hdots&m_1(\rho^{N-1}z)\\
\vdots&\vdots&\vdots&\vdots\\
m_{N-1}(z)&m_{N-1}(\rho z)&\hdots&m_{N-1}(\rho^{N-1}z)
\end{array}\right) \mbox{ is unitary for a.e }z\in\mathbb{T},
\end{equation}
($\rho=e^{\frac{2\pi i}{N}}$), and define $\psi_i\in H$ by
\begin{equation}\label{eq1_3}
\psi_i=:U^{-1}\pi(m_i)\varphi,\quad(i\in\{1,...,N-1\}).
\end{equation}
Then
\begin{equation}\label{eq1_4}
\{T^k\psi_i\,|\,k\in\mathbb{Z},i\in\{1,...,N-1\}\}\mbox{ is an
orthonormal basis for }V_1\ominus V_0
\end{equation}
and
\begin{equation}\label{eq1_5}
\{U^mT^n\psi_i\,|\,m,n\in\mathbb{Z},i\in\{1,...,N-1\}\}\mbox{ is
an orthonormal basis for }H.
\end{equation}
\end{proposition}

\par
Consider a Lipschitz function $m_0$ on $\mathbb{T}$ that satisfies
the following conditions
\begin{equation}\label{eq2_1}
R_{m_0}1=1,
\end{equation}
\begin{equation}\label{eq2_2}
m_0\mbox{ has a finite number of zeros}.
\end{equation}
We call a cycle $C=\{z_1=e^{-i\theta_1},...,z_p=e^{-i\theta_p}\}$
an {\it $m_0$-cycle} if $|m_0(z_k)|=\sqrt{N}$ for all
$k\in\{1,...,p\}$.
\par
We assume in addition that
\begin{equation}\label{eq2_3}
\mbox{ There is at least one }m_0\mbox{-cycle}.
\end{equation}
\par
We have shown in \cite{Dut2}, Proposition 2.13 that for each
$m_0$-cycle one can construct a scaling vector with filter $m_0$ in
the wavelet representation $\mathfrak{R}_{C,\alpha}$ where
$\alpha_k=m_0(z_k)/\sqrt{N}$, $(k\in\{1,...,p\})$. The scaling
vector is defined in the Fourier space as an infinite product:
\begin{equation}\label{eq2_4}
\widehat{\varphi}_{C,k}(x):=\prod_{l=1}^\infty\frac{\cj\alpha_{k-l}m_0\left(\frac{x}{N^l}+\theta_{k-l}\right)}{\sqrt{N}},\,(x\in\mathbb{R}),\varphi_C:=(\varphi_{C,1},...,\varphi_{C,p}),
\end{equation}
(here the subscripts of $\theta$ are considered modulo $p$, that is
$\theta_0=\theta_p,\theta_{-1}=\theta_{p-1},z_{p+2}=z_2$, etc.).
When $m_0$ is a trigonometric polynomial each component of the
scaling vector is compactly supported (see the argument used in
\cite{Dau92}, lemma 6.2.2).
\par
As explained in \cite{Dut2}, Proposition 2.13, the function
$$h_C(\theta):=h_{\varphi_C}(\theta)=\sum_{k=1}^p\operatorname*{Per}|\widehat{\varphi}_{C,k}|^2(\theta-\theta_k),\quad(\theta\in[-\pi,\pi]),$$
is non-negative, harmonic for $R_{m_0}$ and Lipschitz
(trigonometric polynomial when $m_0$ is one).
\par
Moreover $h_C$ is constant $1$ on the $m_0$-cycle $C$ and it is
constant $0$ on every other $m_0$-cycle. This makes $h_C$ linearly
independent for different $m_0$-cycles.
\par
By remark 5.2.4 in \cite{BraJo}, the dimension of the eigenspace
$$\{h\in C(\mathbb{T})\,|\,R_{m_0}h=h\}$$ is equal to the
number of $m_0$-cycles so that
$$\{h_C\,|\, C\mbox{ is an }m_0\mbox{-cycle}\}$$
is a basis for this eigenspace. This shows that (ii) and (iii) in
the next theorem are equivalent. Using Theorem 2.16 in \cite{Dut2}
we obtain:
\begin{theorem}\label{th2_1}
Let $m_0$ be a Lipschitz function satisfying
(\ref{eq2_1}),(\ref{eq2_2}) and (\ref{eq2_3}) and let $C_1,...,C_n$
be distinct $m_0$-cycles. Denote by $\alpha_{ij}$ the coefficients
$m_0(z_{ij})/\sqrt{N}$ where for each fixed $i$ the numbers $z_{ij}$
run through the cycle $C_i$. Define $\varphi_{C_i}$ as in
(\ref{eq2_4}) $(i\in\{1,...,n\})$. Then
$\varphi_{C,\alpha}:=\varphi_{C_1}\oplus...\oplus\varphi_{C_n}$ is a
scaling vector for the wavelet representation
$\mathfrak{R}_{C,\alpha}:=\mathfrak{R}_{C_1,\alpha_1}\oplus...\oplus\mathfrak{R}_{C_n,\alpha_n}$.
Moreover, the following affirmations are equivalent:
\begin{enumerate}
\item
$\varphi_{C,\alpha}$ is an orthogonal scaling vector;
\item {\bf [Cohen's condition]}
The number of $m_0$-cycles is $n$;
\item {\bf [Lawton's condition]}
The dimension of the eigenspace
$$\{h\in C(\mathbb{T})\,|\,R_{m_0}h=h\}$$
is $n$.
\par
If $n$ is smaller then the number of $m_0$-cycles and
$\psi_1,...,\psi_{N-1}$ are defined as in Proposition
\ref{prop1_11}, then
\begin{equation}\label{eqadd1}
\{U^jT^k\psi_i\,|\,i\in\{1,...,N-1\},j,k\in\mathbb{Z}\},
\end{equation}
is a normalized tight frame for $\ltwor^n$.
\end{enumerate}
\end{theorem}
\begin{proof}
For the last statement we use the following argument: if $n$ is
equal to the number of cycles then $\varphi$ is an orthogonal
scaling function so the family in (\ref{eqadd1}) is an
orthonormal basis. If we project this orthonormal basis onto some
of the components corresponding to a choice of a subset of
$m_0$-cycles, we get the normalized tight frame that corresponds
to the case when $n$ is smaller than the number of cycles.

\end{proof}

\begin{remark}\label{rem2_3}
Theorem \ref{th2_1} generalizes some well known results of A.
Cohen and W. Lawton for the classical wavelet representation
$\mathfrak{R}_0$ (see \cite{Co90}, \cite{Law91a}, \cite{Dau92}).
However much more is true: each normalized tight frame wavelet
obtained in the way described before is in fact a projection of
an orthonormal super-wavelet, the one which resides in the larger
space of the larger wavelet representation
$$\mathfrak{R}_{C_1,\alpha_1}\oplus...\oplus\mathfrak{R}_{C_m,\alpha_m}$$
which takes into consideration {\em all} the $m_0$-cycles.
\par
Of course, it is clear that this projection lies in the commutant
of the larger representation and so all operators are intertwined.
\end{remark}

\section{\label{casc}Convergence of the cascade algorithm}
We saw that the orthogonal scaling vector must satisfy the scaling
equation
$$U\varphi=\pi(m_0)\varphi.$$
\par
Generically, there is no closed formula for the scaling
function/vector. For some choices of the filter $m_0$, the scaling
function can have a fractal nature (see e.g. \cite{BraJo}). In
applications, numerical values are needed. The cascade algorithm
provides approximates of the scaling vector $\varphi$ in the
$L^2$-norm. It starts with a well chosen function $\psi^{(0)}$, and,
by iteration of the refinement (or cascade) operator
$$M:=U^{-1}\pi(m_0),$$
$$\psi^{(n+1)}:=M\psi^{(n)},$$
it produces a sequence which converges towards the scaling vector
$\varphi$:
\begin{equation}\label{eqcasc1}
\lim_{n\rightarrow\infty}\|\varphi-\psi^{(n)}\|=0.
\end{equation}
\par
The question here is: how should $\psi^{(0)}$ be chosen so that
the algorithm is convergent to the scaling vector? We will answer
this question in this section. The result is related to spectral
properties of the transfer operator $R_{m_0}$. For a treatment
of the $\ltwor$ case we refer to \cite{BraJo99}.
\par
We will consider $m_0$, a Lipschitz function on $\mathbb{T}$
satisfying (\ref{eq2_1}), (\ref{eq2_2}) and (\ref{eq2_3}). Let
$C_1,...,C_n$ be all the $m_0$-cycles and define the wavelet
representation $\mathfrak{R}_{C,\alpha}$ and the orthogonal scaling
vector $\varphi_{C,\alpha}$ as is Theorem \ref{th2_1}.
\par
We will use the notation
$\{z_{i1}=e^{-i\theta_{i1}},z_{i2}=e^{-i\theta_{i2}},...,z_{ip_i}=e^{-i\theta_{ip_i}}\}$
for the points in the $m_0$-cycle, and the vectors in the Hilbert
space of this representation are of the form
$$(\xi_{ij})_{ij},\quad\mbox{with
}\xi_{ij}\in\ltwor,i\in\{1,...,n\},j\in\{1,...,p_i\}.$$
\par
Consider a starting vector $\psi^{(0)}$ in $\ltwor_{C,\alpha}$ and
define inductively
$$\psi^{(n+1)}:=U_{C,\alpha}^{-1}\pi_{C,\alpha}(m_0)\psi^{(n)}.$$
The result is
\begin{theorem}\label{thc_1}
If
\begin{equation}\label{eqc1_1}
\operatorname*{Per}|\widehat{\psi}^{(0)}_{ij}|^2(\theta_{ij}-\theta_{i'j'})=\delta_{ii'}\delta_{jj'},
\end{equation}
for all
$i,i'\in\{1,...,n\},j\in\{1,...,p_i\},j'\in\{1,...,p_{i'}\}$, and
\begin{equation}\label{eqc1_2}
\widehat{\psi}^{(0)}_{ij}(2k\pi)=\delta_k,\quad(i\in\{1,...,n\},j\in\{1,...,p_i\},k\in\mathbb{Z}),
\end{equation}
then
\begin{equation}\label{eqc1_3}
\lim_{n\rightarrow\infty}\|\varphi-\psi^{(n)}\|=0.
\end{equation}
A simple choice of such a $\psi^{(0)}$ is
\begin{equation}\label{eqc1_4}
\psi_{ij}=\frac{1}{L}\chi_{[0,L)},\quad(i\in\{1,...,n\},j\in\{1,...,p_i\}),
\end{equation}
where $L$ is a positive integer with the property that
$z_{ij}^L=1$ for all $i,j$.
\end{theorem}
\begin{proof}
The proof of the theorem will involve the use of some spectral
properties of the transfer operator $R_{m_0}$ regarded as an
operator on $C(\mathbb{T})$.
\par
By Theorem 3.4.4. in \cite{BraJo}, $R_{m_0}$ has a finite number of
eigenvalues $\lambda_1,...,\lambda_p$ of modulus 1 and $R_{m_0}$ has
a decomposition
\begin{equation}\label{eqj1}
R_{m_0}=\sum_{i=1}^p\lambda_iT_{\lambda_i}+S,
\end{equation}
where $T_{\lambda_i}$ and $S$ are bounded operators on
$C(\mathbb{T})$ such that
\begin{equation}\label{eqj2}
T_{\lambda_i}^2=T_{\lambda_i},\quad
T_{\lambda_i}T_{\lambda_j}=0,\mbox{ for }i\neq j,\quad
T_{\lambda_i}S=ST_{\lambda_i}=0.
\end{equation}
There is a constant $M>0$ such that
\begin{equation}\label{eqj3}
\|S^n\|\leq M,\quad(n\in\mathbb{N}).
\end{equation}
\par
By Proposition 4.4.4 in \cite{BraJo}, the following conditions
(\ref{eqj4}) and (\ref{eqj5}) are equivalent:
\begin{equation}\label{eqj4}
\lim_{n\rightarrow\infty}R_{m_0}^n(g)=T_{1}(g);
\end{equation}
\begin{equation}\label{eqj5}
T_{\lambda_i}(g)=0,\quad\mbox{for all }\lambda_i\neq 1
\end{equation}
\par
By Theorem 2.17 and Corollary 2.18 in \cite{Dut2}, a $\lambda$ with
$|\lambda|=1$ is an eigenvalue for $R_{m_0}$ if and only if there is
an $i\in\{1,...,n\}$ such that $\lambda^{p_i}=1$.
\par
For such a $\lambda$, the eigenspace is described as follows: for
each $i\in\{1,...,n\}$ with the property that $\lambda^{p_i}=1$,
one can define a continuous function $h_i^{\lambda}$ with
$$R_{m_0}h_i^{\lambda}=\lambda h_i^{\lambda},$$
so that for different indices $i$ the functions $h_i^{\lambda}$
are linearly independent and, in fact, they form a basis for the
eigenspace
$$\{h\in C(\mathbb{T})\,|\, R_{m_0}h=\lambda h\}.$$
\par
Define the discrete measures
\begin{equation}\label{eqj6}
\nu_i^\lambda:=\frac{1}{p_i}\sum_{j=1}^{p_i}\lambda^{j-1}\delta_{z_{ij}},\quad(i\in\{1,...,n\},\lambda\in\mathbb{T},\lambda^{p_i}=1),
\end{equation}
where $\delta_{z}$ is the Dirac measure at $z$. Then
\begin{equation}\label{eqj7}
T_{\lambda}(f)=\sum_{i\in\{1,...,n\}\mbox{ with
}\lambda^{p_i}=1}\nu_i^\lambda(f)h_i^\lambda.
\end{equation}
\par
These are the tools that we need for the proof of Theorem
\ref{thc_1}. \par In the proof we will omit the subscript $C,\alpha$
to simplify the notation. We have the following relation between
successive approximations: for $f\in\linft$
\begin{align*}
\int_{\mathbb{T}}fh_{\varphi-\psi^{(n+1)}}\,d\mu&=\ip{\pi(f)(\varphi-\psi^{(n+1)})}{\varphi-\psi^{(n+1)}}\\
&=\ip{\pi(f)U^{-1}\pi(m_0)(\varphi-\psi^{(n)})}{U^{-1}\pi(m_0)(\varphi-\psi^{(n)})}\\
&=\ip{\pi(f(z^N)m_0)(\varphi-\psi^{(n)})}{\pi(m_0)(\varphi-\psi^{(n)})}\\
&=\int_{\mathbb{T}}f(z^N)|m_0|^2h_{\varphi-\psi^{(n)}}\,d\mu\\
&=\int_{\mathbb{T}}f(z)R_{m_0}h_{\varphi-\psi^{(n)}}\,d\mu.
\end{align*}
Thus, as $f$ is arbitrary:
$$h_{\varphi-\psi^{(n+1)}}=R_{m_0}h_{\varphi-\psi^{(n)}}.$$
This implies by induction that
\begin{align*}
\ip{\varphi-\psi^{(n)}}{\varphi-\psi^{(n)}}&=\int_{\mathbb{T}}h_{\varphi-\psi^{(n)}}\,d\mu\\
&=\int_{\mathbb{T}}R_{m_0}h_{\varphi-\psi^{(n-1)}}\,d\mu=...=\int_{\mathbb{T}}R_{m_0}^n(h_{\varphi-\psi^{(0)}})\,d\mu.
\end{align*}
So
\begin{equation}\label{eq6_12_1}
\|\varphi-\psi^{(n)}\|^2=\int_{\mathbb{T}}R_{m_0}^n(h_{\varphi-\psi^{(0)}})\,d\mu.
\end{equation}
\par
Note that,
$$h_{\varphi-\psi^{(0)}}= h_{\varphi}-2\operatorname*{Re}h_{\varphi,\psi^{(0)}}+ h_{\psi^{(0)}},$$
because, for all real-valued $f\in\linft$,
$$\ip{\pi(f)(\varphi-\psi^{(0)})}{\varphi-\psi^{(0)}}=\ip{\pi(f)\varphi}{\varphi}-2\operatorname*{Re}\ip{\pi(f)\varphi}{\psi^{(0)}}+\ip{\pi(f)\psi^{(0)}}{\psi^{(0)}}.$$
\par
We know that $h_{\varphi}=1$ because $\varphi$ is an orthogonal
scaling vector.
\par
We see that the correlation function
\begin{equation}\label{eqc1_5}
h_{\varphi,\psi^{(0)}}(\theta)=\sum_{i=1}^n\sum_{j=1}^{p_i}\operatorname*{Per}(\widehat{\varphi}_{ij}\cj{\widehat{\psi}^{(0)}_{ij}})(\theta-\theta_{ij}),
\end{equation}
\begin{equation}\label{eqc1_6}
h_{\psi^{(0)}}(\theta)=\sum_{i=1}^n\sum_{j=1}^{p_i}\operatorname*{Per}|\widehat{\psi}^{(0)}_{ij}|^2(\theta-\theta_{ij}).
\end{equation}
\par
We have to compute $T_{\lambda}(h_{\varphi,\psi^{(0)}})$ and
$T_{\lambda}(h_{\psi^{(0)}})$ which amounts to computing the
values of $h_{\varphi,\psi^{(0)}}$ and $h_{\psi^{(0)}}$ on the
$m_0$-cycles. For this we will use the inequality
\begin{equation}\label{eqc1_7}
|\operatorname*{Per}(\widehat{\varphi}_{ij}\cj{\widehat{\psi}^{(0)}_{ij}})(\theta)|^2\leq
\operatorname*{Per}|\widehat{\varphi}_{ij}|^2(\theta)
\operatorname*{Per}|\widehat{\psi}^{(0)}_{ij}|^2(\theta).
\end{equation}
Also, we know (see \cite{Dut2}, Proposition 2.13 and its proof) that
the function
$$g_{ij}(\theta):=\operatorname*{Per}|\widehat{\varphi}_{ij}|^2(\theta-\theta_{ij}),\quad(i\in\{1,...,n\},j\in\{1,...,p_i\}),$$
has
\begin{equation}\label{eqc1_8}
g_{ij}(\theta_{kl})=\delta_{ik}\delta_{jl},
\end{equation}
and also
\begin{equation}\label{eqc1_9}
\widehat{\varphi}_{ij}(2k\pi)=\delta_k\quad(k\in\mathbb{Z}).
\end{equation}
Then, inserting (\ref{eqc1_1}), (\ref{eqc1_2}) in (\ref{eqc1_6}),
we get
\begin{equation}\label{eqc1_10}
h_{\psi^{(0)}}(\theta_{kl})=\operatorname*{Per}|\widehat{\psi}^{(0)}_{kl}|^2(0)=1.
\end{equation}
Using in (\ref{eqc1_5}) the relations (\ref{eqc1_7}),
(\ref{eqc1_8}) and then (\ref{eqc1_2}) and (\ref{eqc1_9}), we have
\begin{equation}\label{eqc1_11}
h_{\varphi,\psi^{(0)}}(\theta_{kl})=
\operatorname*{Per}(\widehat{\varphi}_{kl}\cj{\widehat{\psi}^{(0)}_{kl}})(0)=1.
\end{equation}
With (\ref{eqc1_10}) and (\ref{eqc1_11}) in (\ref{eqj7}) and
(\ref{eqj6}), it follows that
$$T_{\lambda}(h_{\psi^{(0)}})=\delta_{\lambda,1}\cdot 1\quad
T_{\lambda}(\operatorname*{Re}h_{\varphi,\psi^{(0)}})=\delta_{\lambda,1}\cdot
1.$$ Therefore (\ref{eqj5}) is satisfied and
$$\lim_{n\rightarrow\infty}R_{m_0}^n(h_{\psi^{(0)}})=
\lim_{n\rightarrow\infty}R_{m_0}^n(\operatorname*{Re}h_{\varphi,\psi^{(0)}})=1$$
which implies, with (\ref{eq6_12_1}), that
$$\lim_{n\rightarrow\infty}\|\varphi-\psi^{(n)}\|^2=0.$$
\par
It only remains to check that the vector $\psi^{(0)}$ given in
equation (\ref{eqc1_4}) satisfies (\ref{eqc1_1}) and
(\ref{eqc1_2}).
$$\widehat{\psi}_{ij}^{(0)}(\xi)=e^{-iL\xi/2}\frac{\sin(L\xi/2)}{L\xi/2},$$
hence
$$\widehat{\psi}_{ij}^{(0)}(\frac{2k\pi}{L})=\delta_k,\quad(k\in\mathbb{Z}).$$
Since $z_{ij}^L=1$, we see that
$e^{-iL(\theta_{ij}-\theta_{i'j'})}=1$ so
$\theta_{ij}-\theta_{i'j'}=\frac{2k_0\pi}{L}$ for some
$k_0\in\mathbb{Z}$. Therefore, if $(i,j)\neq (i',j')$ then
$k_0\neq 0\,\mod L$ and
$$\operatorname*{Per}|\widehat{\psi}^{(0)}_{ij}|^2(\theta_{ij}-\theta_{i'j'})=
\sum_{k\in\mathbb{Z}}|\widehat{\psi}^{(0)}_{ij}|^2(\frac{2k_0\pi+2kL\pi}{L})=0.
$$
Also
$$\operatorname*{Per}|\widehat{\psi}^{(0)}_{ij}|^2(0)=1$$
and
$$\widehat{\psi}^{(0)}_{ij}(2k\pi)=\delta_k,\quad(k\in\mathbb{Z}).$$
\end{proof}

\section{\label{exam}Examples}

\begin{example}\label{ex1_15}
Consider the wavelet representation $\mathfrak{R}_0$ and suppose
$m_0\in\linft$ is a filter that has an orthogonal scaling function
$\varphi\in\ltwor$ for this representation.
\par
Define a new filter $\tilde m_0\in\linft$ as follows: let $p$ be a
positive integer which is prime with $N$, and define
$$\tilde{m_0}(z)=m_0(z^p),\quad(z\in \mathbb{T}).$$
Because $p$ and $N$ are mutually prime, the $p$-th roots of unity,
$\{z\in\mathbb{T}\,|\,z^p=1\}$ split into several disjoint cycles
(the map $z\mapsto z^N$ is a bijection on
$\{z\in\mathbb{T}\,|\,z^p=1\}$); for example, if $N=2$, $p=9$ and
$\rho_k=e^{-i\frac{2k\pi}{p}}$, $i\in\{0,...,8\}$ are the $9$-th
roots of 1, then the cycles are
$$\{\rho_0\},\quad\{\rho_1,\rho_2,\rho_4,\rho_8,\rho_7,\rho_5\},\quad\{\rho_3,\rho_6\}.$$
\par
Let $C_1,...,C_n$ be these cycles
$$C_1\cup...\cup C_n=\{z\in\mathbb{T}\,|\, z^p=1\}.$$
We show that $\tilde m_0$ has an orthogonal scaling vector for the
wavelet representation
$$\mathfrak{R}_{C_1}\oplus...\oplus\mathfrak{R}_{C_n},$$
namely
$$\tilde\varphi:=(\tilde\varphi_0,...,\tilde\varphi_p),\quad
\tilde\varphi_i(x)=\frac{1}{p}\varphi\left(\frac{x}{p}\right),\quad(x\in\mathbb{R},i\in\{0,...,p-1\}).$$
\par
We have to verify the conditions of Corollary \ref{cor1_15}. For
this we first write the conditions that are satisfied by the scaling
function $\varphi$:
\begin{equation}\label{eqex1_15_1}
\operatorname*{Per}|\widehat{\varphi}|^2(x)=1,\quad(x\in\mathbb{R});
\end{equation}
\begin{equation}\label{eqex1_15_2}
\sqrt{N}\widehat{\varphi}(Nx)=m_0(x)\widehat{\varphi}(x),\quad(x\in\mathbb{R});
\end{equation}
\begin{equation}\label{eqex1_15_3}
\mbox{There is no }N\mbox{-invariant set of positive measure such
that }\widehat{\varphi}\mbox{ vanishes on it.}
\end{equation}
Having these, we check the conditions for $\tilde\varphi$.
\par
The orthogonality condition can be restated in this case as
$$
\sum_{j=0}^{p-1}\operatorname*{Per}|\widehat{\tilde{\varphi}}|^2(x-\frac{2j\pi}{p})=1\quad
(x\in\mathbb{R}).$$ This holds because \[
\sum_{j=0}^{p-1}\operatorname*{Per}|\widehat{\tilde{\varphi}}|^2(x-\frac{2j\pi}{p})=\sum_{i=0}^{p-1}
\sum_{k\in \mathbb{Z}}|\widehat{\varphi}|^2(px-2j\pi+2k\pi p)=\]
\[\sum_{k\in\mathbb{Z}}|\widehat{\varphi}|^2(px+2k\pi)=\operatorname*{Per}|\widehat{\varphi}|^2(px)=1\quad
(x\in\mathbb{R}),\] where we used (\ref{eqex1_15_1}) in the last
equality.
\par
For the scaling equation we have to check that
$$\sqrt{N}\widehat{\varphi}(pNx)=m_0\left(p\left(\frac{2j\pi}{p}+x\right)\right)\widehat\varphi(px),\quad(x\in\mathbb{R},j\in\{0,...,p-1\}),$$
which is clear from (\ref{eqex1_15_2}).
\par
Suppose now that the cyclicity condition is not satisfied. Then some
of the components of $\tilde\varphi$ that correspond to one of the
cycles do not satisfy the cyclicity condition. Let $l$ be the length
of this cycle. Then there is an $N^l$-invariant set of positive
measure, call it $E$, such that $\widehat{\varphi}(px)$ vanishes on
$E$.
\par
Since $\tilde\varphi$ satisfies the scaling equation, it follows
that $\widehat\varphi(px)$ vanishes also on
$NE,N^2E,...,N^{l-1}E$.
\par
Take
$$A=\frac{1}{p}(E\cup NE\cup...\cup N^{l-1}E).$$
Then $A$ is $N$-invariant, of positive measure, and
$\widehat\varphi(x)$ vanishes on $A$. This contradicts
(\ref{eqex1_15_3}).
\par
In conclusion, $\tilde\varphi$ is indeed an orthogonal scaling
vector with filter $\tilde m_0$ for the wavelet representation
$$\mathfrak{R}_{C_1}\oplus...\oplus\mathfrak{R}_{C_n}.$$
\end{example}

\par
The next example shows that for the scale $N=2$, no matter how we
choose the cycles $C=C_1\cup C_2\cup...\cup C_p$, there exists a
MRA, orthogonal super-wavelet for the representation
$\mathfrak{R}_C$.

\begin{example}\label{ex6_2} Let $C_1,C_2,...,C_p$ be 2-cycles and let $C=C_1\cup C_2\cup ...\cup C_p$. We will
construct an orthogonal scaling vector $\varphi$ for the wavelet
representation
$$\mathfrak{R}_{C_1}\oplus...\oplus\mathfrak{R}_{C_p}.$$
The following definitions for $x\in \mathbb{R}$
will be used in this example :

\begin{enumerate}
\item $x$ is called a cycle point if there is $c$ in
$C$ such that $$x\equiv \theta\mbox{ mod }2\pi, \mbox{ where
}e^{-i\theta}=c;$$
\item $x$ is called a supplement if $x-\pi$ is a cycle point;
\item $x$ is called a main point if it is a cycle point or a
supplement;
\item $x$ is called  mid-point if $x=\frac{a+b}{2}$, with $a,b$
consecutive main points;
\item $x$ is called a cycle midpoint if $x=\frac{a+b}{2}$ with
$a,b$ consecutive cycle points.
\end{enumerate}
(Here, when we say ``consecutive", we refer to the order on the real
line ).
\par
For $z\in C, z=e^{-i\theta_0}$, $\theta\in[-\pi,\pi]$, define
$$\widehat{\varphi_z}(\theta)=\chi_{[\frac{a(\theta_0)+\theta_0}{2},\frac{\theta_0+b(\theta_0)}{2}]}(\theta+\theta_0),$$
where $a(\theta_0),\theta_0,b(\theta_0)$ are consecutive cycle
points. It is easy to check that
$$\sum_{z=e^{-i\theta_0}\in C}\operatorname*{Per}|\widehat{\varphi}_z|^2(\theta-\theta_0)=1,\quad(\theta\in\mathbb{R}).$$
Hence $\varphi$ defined by its Fourier transform
$$\widehat{\varphi}=\oplus_{z\in
C}\widehat{\varphi}_z=\oplus_{i=1}^p\oplus_{z\in
C_i}\widehat{\varphi}_z=:\oplus_{i=1}^p\widehat{\varphi}_{C_i}$$
is a good candidate for an orthogonal scaling vector corresponding
to $C$ build out of orthogonal scaling vectors corresponding to
each $C_i$ . Next let us define the filter $m_0$ :
\[m_0=\sum_{\theta_0\in [-\pi,\pi], cycle point
}\chi_{\left([\frac{c(\theta_0)+\theta_0}{2},\frac{\theta_0+d(\theta_0)}{2}]\cap[-\pi,\pi]\right)},
\]
where for the cycle point $\theta_0\in [-\pi,\pi],
c(\theta_0),\theta_0,d(\theta_0)$ are consecutive main points. We
will now check the scaling equation for the above defined
$\varphi$ and filter $m_0$. It suffices to show that for two
consecutive elements $z_0=e^{i\theta_0},z_1=e^{i\theta_1}$ of a
cycle $C_i$ (i.e. $z_0^2=z_1$, or equivalently
$2\theta_0\equiv\theta_1$ mod $2\pi$) the following holds :
$$\widehat{\varphi}_{z_1}(2\theta)=m_0(\theta+\theta_0)\widehat{\varphi}_{z_0}(\theta),
\mbox{ a.e. }\theta \in \mathbb{R}$$ or, equivalently
\begin{equation}\label{eq6_2_1}
\widehat{\varphi}_{z_1}(2\theta-2\theta_0)=m_0(\theta)\widehat{\varphi}_{z_0}(\theta-\theta_0),
\mbox{ a.e. }\theta \in \mathbb{R}.
\end{equation}
 Suppose
$\widehat{\varphi}_{z_0}(\theta-\theta_0)=\chi_{[\frac{\alpha
+\theta_0}{2},\frac{\beta + \theta_0}{2}]}(\theta)$, where
$\alpha<\theta_0<\beta$ are consecutive cycle points. By the
definition of $m_0$, it follows that there are consecutive main
points $a<\theta_0<b$ such that
\begin{equation}\label{eq6_2_2}
m_0(\theta)\widehat{\varphi}_{z_0}(\theta-\theta_0)=\chi_{[\frac{a
+\theta_0}{2},\frac{b + \theta_0}{2}]}(\theta).
\end{equation}
Actually $a$ is either $\alpha$ or the first supplement on the
left of $\theta_0$ and $b$ is either $\beta$ or the first
supplement on the right of $\theta_0$. Suppose
$\widehat{\varphi}_{z_1}(\theta-\theta_1)=\chi_{[\frac{a_1
+\theta_1}{2},\frac{b_1 + \theta_1}{2}]}(\theta)$, with
$a_1<\theta_1<b_1$ consecutive cycle points. Then\[
\widehat{\varphi}_{z_1}(\theta)=\chi_{[\frac{a_1
-\theta_1}{2},\frac{b_1 - \theta_1}{2}]}(\theta)\] and
\begin{equation}\label{eq6_2_3}
\widehat{\varphi}_{z_1}(2\theta-2\theta_0)=\chi_{[\frac{a_1
-\theta_1}{2},\frac{b_1 - \theta_1}{2}]}(2\theta-2\theta_0).
\end{equation}
Note that, under the map $x\mapsto 2x$, the consecutive main points
$a<\theta_0<b$ are mapped into consecutive cycle points. Since
$2\theta_0 \equiv \theta_1$ mod $2\pi$, and because a translation by
an integer multiple of $2\pi$ maps consecutive cycle points to
consecutive cycle points, it follows in particular that there exists
$k\in \mathbb{Z}$ such that $a_1=2a+2k\pi$,
$\theta_1=2\theta_0+2k\pi$ and $b_1=2b+2k\pi$. We have:
\[\frac{a_1-\theta_1}{2}\le2\theta-2\theta_0\le
\frac{b_1-\theta_1}{2}\Leftrightarrow\]
\[
\frac{2a+2k\pi-(2\theta_0+2k\pi)}{2}\le2\theta-2\theta_0\le\frac{2b+2k\pi-(2\theta_0+2k\pi)}{2}\Leftrightarrow\]
\[\frac{a+\theta_0}{2}\le\theta\le\frac{b+\theta_0}{2},\]
and, with (\ref{eq6_2_2}), (\ref{eq6_2_3}), the relation
(\ref{eq6_2_1}) is obtained.
\par
The cyclicity condition is automatically satisfied because all
$\widehat{\varphi}_z$ contain a neighborhood of $0$.
\par
Consequently $\varphi$ is an orthogonal scaling vector with filter
$m_0$.
\end{example}
We can use Example \ref{ex6_2} to obtain some information about the
distribution of cycles. The idea is that when $m_0$ is the
characteristic function of some intervals, we can obtain a Cohen
condition as in Theorem \ref{th2_1}.
\begin{proposition}\label{prop6_3}
Suppose $m_0$ is of the form
$$m_0=\sqrt{N}\chi_E,$$
where $E\subset[-\pi,\pi]$ is a union of intervals such that none of
the endpoints of these intervals lies on a cycle. Assume moreover
that for some distinct cycles $C_1,...,C_n$, the wavelet
representation
$$\mathfrak{R}_C:=\mathfrak{R}_{C_1}\oplus...\oplus\mathfrak{R}_{C_n}$$
has an orthogonal scaling vector with filter $m_0$.
\par
Then every cycle $D$ disjoint from $C_1\cup...\cup C_n$ must have a
point in $[-\pi,\pi]\setminus E$.
\end{proposition}

\begin{proof}
Since there is an orthogonal scaling vector with filter $m_0$, the
condition
$$R_{m_0}1=1$$
must be satisfied.
\par
Suppose there is a cycle
$D=\{e^{-i\theta_1},...,e^{-i\theta_p}\}$, different from the
given ones such that $D$ is contained in $E$. Because the endpoint
of the intervals of $E$ are not on cycles, each point of $D$ lies
in the interior of $E$. Then one can construct a scaling vector
with filter $m_0$ for the wavelet representation $\mathfrak{R}_D$,
 $\widehat{\varphi}_{D,k}$ defined as in (\ref{eq2_4}). The
fact that $\widehat{\varphi}_{D,k}$ is a well defined $\ltwor$
function can be proved using the arguments in \cite{Dut2},
Proposition 2.13 and \cite{Dau92}, lemma 6.2.1.
\par
The scaling equation can be verified instantly and, since each point
of $D$ lies in the interior of one of the intervals of $m_0$, it
follows that $\widehat{\varphi}_{D,k}$ is $1$ in a neighborhood of
$0$ so the cyclicity condition of Corollary \ref{cor1_15} is also
clear.
\par
But then, having a scaling vector with filter $m_0$ in the wavelet
representation $\mathfrak{R}_D$ means that $\mathfrak{R}_D$ is the
wavelet representation associated to the harmonic function
$h_{\widehat{\varphi}_D}$, the correlation function of
$\varphi_D$ (\cite[Theorem 2.4]{Jor01}).
\par
$h_{\widehat{\varphi}_{D,k}}$ must be bounded for the following
reason: clearly $\widehat{\varphi}_{D,k}$ is either $1$ or $0$.
Also, it is impossible to have an $x\in\mathbb{R}$ with
$\widehat{\varphi}_{D,k}(x)=1$ and
$\widehat{\varphi}_{D,k}(2l_0\pi)=1$ for some integer $l_0\neq0$.
Indeed, otherwise we would have (from (\ref{eq2_4})):
$$m_0\left(\frac{x}{N^l}+\theta_{k-l}\right)=\sqrt{N},\quad
m_0\left(\frac{x+2l_0\pi}{N^l}+\theta_{k-l}\right)=\sqrt{N},\quad(l\in\mathbb{Z},l\geq
1).$$ Write $l_0=N^rq$ with $r,q\in\mathbb{N}$, $q$ not divisible
by $N$. Then
$$m_0\left(\frac{x}{N^{r+1}}+\theta_{k-l}\right)=\sqrt{N},
m_0\left(\frac{x}{N^{r+1}}+\theta_{k-l}+\frac{2l_0\pi}{N}\right)=\sqrt{N},$$
which contradicts $R_{m_0}1=1$.
\par
Thus, if $\widehat{\varphi}_{D,k}(x)=1$ then
$\widehat{\varphi}_{D,k}(x+2k\pi)=0$ for all integers $k\neq0$.
This implies that
$\operatorname*{Per}|\widehat{\varphi}_{D,k}|^2(x)\leq 1$ so
$$h_{\widehat{\varphi}_{D,k}}(x)=\sum\operatorname*{Per}|\widehat{\varphi}_{D,k}|^2(x-\theta_k)\leq
p.$$
\par
Having these, with \cite[Theorem 2.4]{Dut1}, there exists some positive operator $S$ that commutes with $U_C$ and $\pi_C$ such that
the correlation function $h_{S\varphi_C,\varphi_C}$ is $h_{\widehat{\varphi}_{D,k}}$. Then the correlation function of the vector $S^{1/2}\varphi_C$ is
$h_{\widehat{\varphi}_{D,k}}$. Mapping $U_D^{-n}\pi_D(f)\varphi_D$ to $U_C^{-n}\pi_C(f)S^{1/2}\varphi_C$, we get a non-trivial operator which intertwines the representations associated to $D$ and $C$.
But this is impossible because the representations are disjoint \cite[Lemma 2.14]{Dut2}.
\par
In conclusion the assumption was erroneous so the cycle $C$ must
intersect the complement of $E$.
\end{proof}

We illustrate now some particular cases of the filters, orthogonal
scaling vectors, and wavelets given in Example \ref{ex6_2}. Here the
scale is $N=2$. Since the low-pass filter $m_0$ is just a
characteristic function, the corresponding high-pass filter $m_1$
can be chosen as the characteristic function of the complement of
the set that gives $m_0$, and (\ref{eq1_2}) holds. Then the wavelet
is defined as in (\ref{eq1_3}).

\begin{example}\label{ex6_2_2}
Consider the cycle $$C=\{e^{-\frac{2\pi i}{7}},e^{-\frac{4\pi
i}{7}},e^{-\frac{-6\pi i}{7}}\}.$$
$$m_0=\sqrt{2}\chi_{[-\pi,-\frac{11\pi}{14}]\cup[\frac{3\pi}{14},\pi]},$$
$$\widehat{\varphi}_1=\chi_{[-\frac{4\pi}{7},\frac{\pi}{7}]},\quad
\widehat{\varphi}_2=\chi_{[-\frac{\pi}{7},\frac{2\pi}{7}]},\quad
\widehat{\varphi}_3=\chi_{[-\frac{2\pi}{7},\frac{4\pi}{7}]}$$ Then
$(\varphi_1,\varphi_2,\varphi_3)$ is an orthogonal scaling vector
for the wavelet representation $\mathfrak{R}_C$ with filter $m_0$.
\par
The high-pass filter is
$$m_1=\sqrt{2}\chi_{[-\frac{11\pi}{14},\frac{3\pi}{14}]},$$
and the orthogonal wavelet is $(\psi_1,\psi_2,\psi_3)$ with
$$\widehat{\psi}_1=\chi_{[\frac{\pi}{7},\frac{8\pi}{7}]},\quad\widehat{\psi}_2=\chi_{[-\frac{8\pi}{7},-\frac{\pi}{7}]},\quad\widehat{\psi}_3=0.$$

\end{example}

\begin{example}\label{ex6_2_8}
Consider the cycles $$C_1=\{e^{-\frac{2\pi i}{5}},e^{-\frac{4\pi
i}{5}},e^{-\frac{-2\pi i}{5}},e^{-\frac{-4\pi i}{5}}\},\quad
C_2=\{e^{-\frac{2\pi i}{3}},e^{-\frac{-2\pi i}{3}}\}.$$
$$m_0=\sqrt{2}\chi_E,$$
where
$$E:=[-\pi,-\frac{19\pi}{30}]\cup[-\frac{15\pi}{30},-\frac{11\pi}{30}]\cup[\frac{11\pi}{30},\frac{15\pi}{30}]\cup[\frac{19\pi}{30},\pi].$$
$$\widehat{\varphi}_1^1=\chi_{[-\frac{6\pi}{15},\frac{2\pi}{15}]},\quad
\widehat{\varphi}_2^1=\chi_{[-\frac{\pi}{15},\frac{3\pi}{15}]},\quad
\widehat{\varphi}_3^1=\chi_{[-\frac{2\pi}{15},\frac{6\pi}{15}]}\quad
\widehat{\varphi}_4^1=\chi_{[-\frac{3\pi}{15},\frac{\pi}{15}]},$$
$$\widehat{\varphi}_1^2=\chi_{[-\frac{2\pi}{15},\frac{\pi}{15}]},\quad
\widehat{\varphi}_2^2=\chi_{[-\frac{\pi}{15},\frac{2\pi}{15}]}.$$
Then $((\varphi_1^1,...,\varphi_4^1),(\varphi_1^2,\varphi_2^2))$
is an orthogonal scaling vector for the wavelet representation
$\mathfrak{R}_{C_1}\oplus\mathfrak{R}_{C_2}$ with filter $m_0$.
\par
The high-pass filter is
$$m_1=\sqrt{2}\chi_{E_1},$$
where
$$E_1:=[-\frac{19\pi}{30},-\frac{15\pi}{30}]\cup[-\frac{11\pi}{30},\frac{11\pi}{30}]\cup[\frac{15\pi}{30},\frac{19\pi}{30}]$$
and the orthogonal wavelet is
$((\psi_1^1,...,\psi_4^1),(\psi_1^2,\psi_2^2))$ with
$$\widehat{\psi}_1^1=\widehat{\psi}_3^1=0,\quad
\widehat{\psi}_2^1=\chi_{[-\frac{12\pi}{15},-\frac{\pi}{15}]\cup[\frac{3\pi}{15},\frac{4\pi}{15}]},\quad
\widehat{\psi}_4^1=\chi_{[-\frac{4\pi}{15},-\frac{3\pi}{15}]\cup[\frac{\pi}{15},\frac{12\pi}{15}]},$$
$$\widehat{\psi}_1^2=\chi_{[\frac{\pi}{15},\frac{4\pi}{15}]},\quad
\widehat{\psi}_2^2=\chi_{[-\frac{4\pi}{15},-\frac{\pi}{15}]}.$$
\end{example}

\begin{acknowledgements}
We would like to thank professor Palle Jorgensen for his suggestions and
constant support.
\end{acknowledgements}


\begin{thebibliography}{abcdef}
\bibitem[BCMO]{BCMO}
L. Baggett, A. Carrey, W. Moran, P. Ohring, {\em General existence
theorems for orthonormal wavelets, an abstract approach} Publ.
Res. Inst. Math. Sci. Kyoto Univ, {\bf 31} (1995), 95-111
\bibitem[BraJo97]{BraJo97}
O. Bratteli, P.E.T. Jorgensen, {\em {I}sometries, shifts, Cuntz
algebras and multiresolution wavelet analysis of scale $N$},
Integral Equations Operator Theory, 28 (1997), 382-443
\bibitem[BraJo99]{BraJo99}
O. Bratteli, P.E.T. Jorgensen, {\em Convergence of the cascade
algorithm at irregular scaling functions}, The Functional and
Harmonic Analysis of Wavelets and Frames (San Antonio, 1999)
(L.W.Baggett and D.R.Larson eds.) Contemp.Math. vol. 247, AMS,
Providence 1999, pp. 93-130
\bibitem[BraJo]{BraJo}
O. Bratteli, P.E.T. Jorgensen, {\em {W}avelets Through a Looking
Glass}, Birkhauser, 2002
\bibitem[Co90]{Co90}
A. Cohen {\em {O}ndelettes , analyses multiresolutions et
traitement numerique du signal } , Ph.D. Thesis , Universite
Paris, Dauphine
\bibitem[Dau92]{Dau92}
I.~Daubechies, {\em {T}en {L}ectures on
{W}avelets}, CBMS-NSF Regional Conf.
  Ser. in Appl. Math., vol.~61, Society for Industrial and Applied Mathematics,
  Philadelphia, 1992.
\bibitem[Dut1]{Dut1}
D.E. Dutkay, {\em {H}armonic analysis of signed Ruelle transfer
operators}, J. Math. Anal. Appl. 273 (2002) 590-617
\bibitem[Dut2]{Dut2}
D.E. Dutkay, {\em {T}he wavelet Galerkin operator}, Journal of
Operator Theory, 51 (2004) 49-70.
\bibitem[Dut3]{Dut3}
D.E. Dutkay, {\em Positive definite maps, representations and
frames}, Reviews in Mathematical Physiscs, vol. 16, No. 4 (2004)
1-27.
\bibitem[Jor01]{Jor01}
P.E.T ~Jorgensen, {\em {R}uelle operators : Functions which are
harmonic with respect to  a transfer operator},
 Mem. Amer. Math. Soc., {\bf 152}, no. 720
\bibitem[HL]{HL}
D. Han, D. Larson, {\em {F}rames, bases and group
representations}, Memoirs of the AMS, sept. 2000, vol.147, no. 697
\bibitem[HeWe]{HeWe}
E. Hernandez, G. Weiss, {\em {A} First Course on Wavelets}, CRC
Press, Inc. 1996
\bibitem[Law91a]{Law91a}
W.M. Lawton, {\em{N}eccesary and sufficient conditions for
constructing orthonormal wavelet bases} , J. Math. Phys. 32
(1991) , 57-61
\end{thebibliography}
\end{document}